\documentclass[final]{amsart}
\usepackage{graphicx}
\usepackage{amssymb}
\usepackage{amsmath}
\usepackage[cmtip,arrow]{xy}
\usepackage{pb-diagram,pb-xy}

\allowdisplaybreaks

\newtheorem{theorem}{Theorem}
\newtheorem{corollary}[theorem]{Corollary}
\newtheorem{lemma}[theorem]{Lemma}

\newtheorem{definition}[theorem]{Definition}

\newtheorem{remark}[theorem]{Remark}

\newtheorem{example}[theorem]{Example}

\begin{document}
%

\title[Stopping Markov processes]{Compressing oversized information
in Markov chains}%
\author{Giacomo Aletti}%
\address{Dipartimento di Matematica \\
Universit\`a\ di Milano, 20133 Milan, ITALY}%
\email{giacomo.aletti\@unimi.it}%


\thanks{The first author wish to thank Ely Merzbach, coauthor of \cite{AM}, who
initially gave him the problem in \cite[Example~2]{AM} which stimulates the present work
and \cite{AM}. He points out one typo, of his responsability. In that example, $25$ should be read $26$.}%

\begin{abstract}
Given a strongly stationary Markov chain
and a finite set of stopping rules, we prove the existence of a
polynomial algorithm which projects the Markov chain onto a
minimal Markov chain without redundant information. Markov complexity
is hence defined and tested on some classical problems.
\end{abstract}

\subjclass{Primary: 60J22; Secondary: 90C35, 94C15}%
\keywords{Markov time of the first passage, stopping rules, Markov complexity}%


%
\date{\today}
%
\maketitle

\section{Introduction}
Let $X_n$ be a stationary Markov chain on a finite set $E$ with
transition matrix $P$. The Markov process stops when one of the
given stopping rules occurs. The problem of finding the stopping
law may be solved by embedding the Markov chain into another
Markov chain on a larger state's set (the tree made by both the
states and the stopping rules, see \cite{AM}). The desired law is
obtained then from the transition matrix of the new Markov chain.

Unfortunately, this new Markov chain may be so big that numerical
computations can be not practicable. A new method permitting to
ensure the existence of a projection of the Markov chain into a
``minimal'' Markov chain which preserves probabilities was
presented in \cite{AM}.

\begin{figure}[tbh]
\begin{center}
  \includegraphics[width=.45\linewidth]{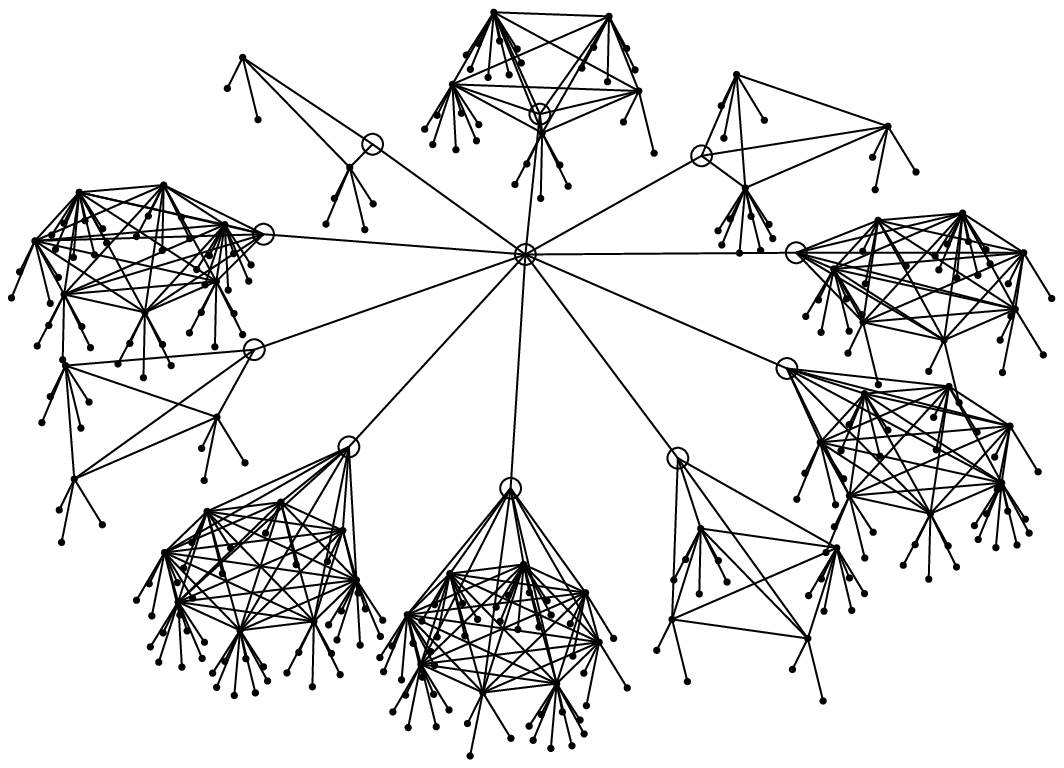}
  \includegraphics[width=.45\linewidth]{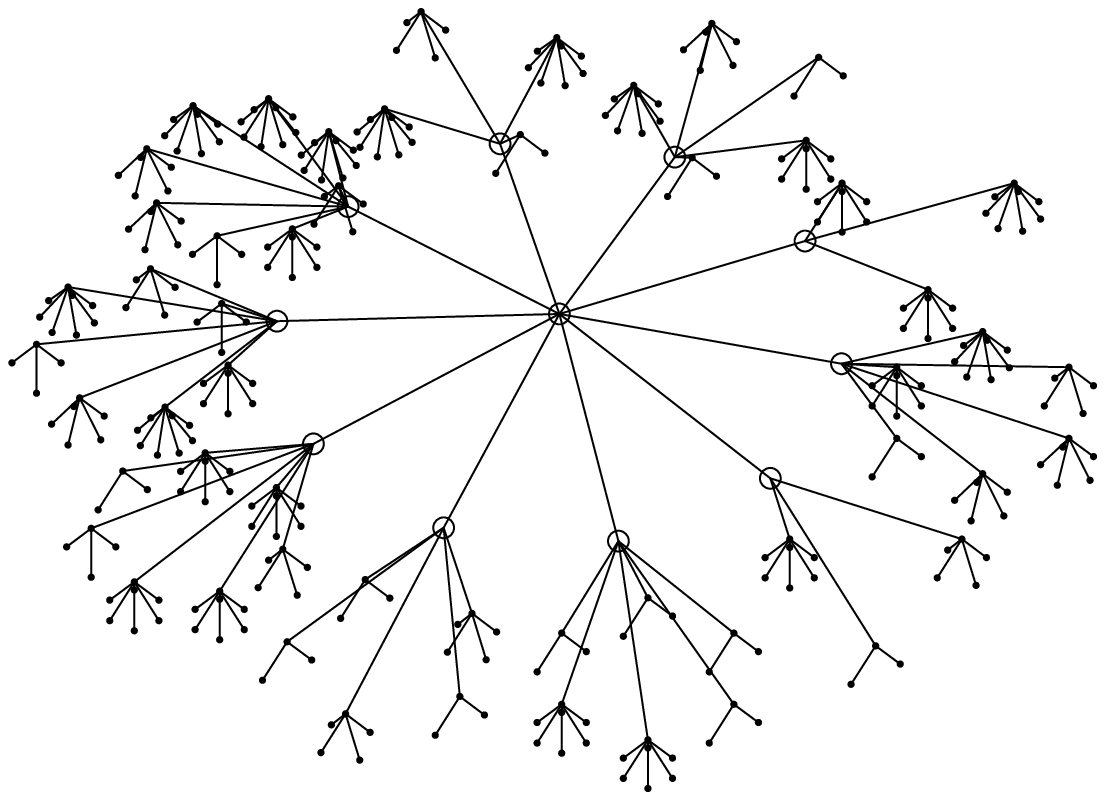}
\end{center}
  \caption{Networks that may be shrinked.}\label{fig:realnet}
\end{figure}

As in \cite{AM}, we recall now how this problem occurs in many
situations.
\begin{enumerate}
    \item In finance some filter rules for trading is a special
case of the Markov chain stopping rule suggested by the authors in
\cite{AM}.
    \item ``When enough is enough''!  For example, an insured has an
accident only occasionally in a while.  How many accidents in a
specified number of years should be used as a stopping time for
the insured (in other words, when it should be discontinued the
insurance contract).
    \item \emph{State dependent markov chains.}
    Namely, the transition probabilities are given in terms of the
history. For simplicity consider the decision to stop if we get
$2$ identical throws ($11$, $22$, $33$, \ldots $nn$) (for example,
when $n=2$, an insured has two kinds of accidents in row-one each
year and is discontinued or an insured has no accidents two years
in a row and therefore he is ``promoted'' to a better class of
insured). If probability of a switch from $hm$ to $mk$ is denoted
by $p_{hm,mk}$ then the Markov transition matrix has the form:
\begin{footnotesize}
\begin{equation}\label{eq:tbcll}
\begin{array}{c|c@{\,}c@{\,}c@{\,}c@{\,}c@{\,}c@{\,}c@{\,}c@{\,}c@{\,}c@{\,}c@{\,}c@{\,}c}
   & 11 & 12 & \ldots & 1n & 21 & 22 & \ldots & 2n & \ldots & n1 & n2 & \ldots & nn \\
   \hline \\[-.2cm]
  11 & p_{11,11} & p_{11,12} & \ldots & p_{11,1n} & 0 & \ldots & 0 & \ldots & 0 & 0 & \ldots & 0 \\
  12 & 0 & 0 & \ldots & 0 & p_{12,21} & p_{12,22} & \ldots & p_{12,2n} & \ldots & 0 & 0 & \ldots & 0 \\
\vdots & \vdots & \vdots & \vdots & \vdots & \vdots & \vdots &
\vdots & \vdots & \vdots & \vdots & \vdots &
\vdots & \vdots \\
  1n & 0 & 0 & \ldots & 0 & 0 & 0 & \ldots & 0 & \ldots & p_{1n,n1} & p_{1n,n2} & \ldots & p_{1n,nn} \\
  21 & p_{21,11} & p_{21,12} & \ldots & p_{21,1n} & 0 & 0 & \ldots & 0 & \ldots & 0 & 0 & \ldots & 0 \\
  22 & 0 & 0 & \ldots & 0 & p_{22,21} & p_{22,22} & \ldots & p_{22,2n} & \ldots & 0 & 0 & \ldots & 0 \\
\vdots & \vdots & \vdots & \vdots & \vdots & \vdots & \vdots &
\vdots & \vdots & \vdots & \vdots & \vdots &
\vdots & \vdots \\
  2n & 0 & 0 & \ldots & 0 & 0 & 0 & \ldots & 0 & \ldots & p_{2n,n1} & p_{2n,n2} & \ldots & p_{2n,nn} \\
\vdots & \vdots & \vdots & \vdots & \vdots & \vdots & \vdots &
\vdots & \vdots & \vdots & \vdots & \vdots &
\vdots & \vdots \\
  n1 & p_{n1,11} & p_{n1,12} & \ldots & p_{n1,1n} & 0 & 0 & \ldots & 0 & \ldots & 0 & 0 & \ldots & 0 \\
  n2 & 0 & 0 & \ldots & 0 & p_{n2,21} & p_{n2,21} & \ldots & p_{n2,2n} & \ldots & 0 & 0 & \ldots & 0 \\
\vdots & \vdots & \vdots & \vdots & \vdots & \vdots & \vdots &
\vdots & \vdots & \vdots & \vdots & \vdots &
\vdots & \vdots \\
  nn & 0 & 0 & \ldots & 0 & 0 & 0 & \ldots & 0 & \ldots & p_{nn,n1} & p_{nn,n2} & \ldots & p_{nn,nn}
\end{array}
\end{equation}
\end{footnotesize}%
which can be analyzed for the stopping time by the usual methods.
Obviously, in many situations (e.g., if $p_{hm,mk}=p_{m,k}$
$\forall h\neq m$), this matrix has a special structure and can be
reduced.
    \item \emph{Small-world Networks}. Given one of the networks as in
    Figure~\ref{fig:realnet}, is it possible to reduce it and to
    preserve the law of reaching a given absorbing state?
\end{enumerate}

Formally, the problem is given by a triple $(E,T,P)$, where:
\begin{itemize}
    \item $E$ is a set (set of states);
    \item $T$ is a nonempty subset of $E$ (target set);
    \item $P:E\times E\to \mathbb R_+$ with the following
    properties:
    \begin{itemize}
        \item $\forall e\in E$, $\sum P(e,\cdot)=1$;
        \item $ P^{-1}(0,\infty) \cap (T\times E)\subseteq (T\times T)$.
    \end{itemize}
    $P$ may be identified to the probability transition matrix $P_{i,j}=P((e_i,e_j))$.
\end{itemize}

As shown in \cite{AM}, compressing non influent information is
equivalent to find a triple $(F,t,P^*)$ and a map $\pi:E\to F$
s.t.
\begin{itemize}
    \item $\pi$ is a surjective set function from $E$ to $F$;
    \item $t=\pi(T)$, $T=\pi^{-1}(t)$;
    \item the following diagram commutes:
\begin{equation}\label{eq:diagram}
\begin{diagram}
   \node{E\times \mathfrak P(E)} \arrow{seee,l}{
   \mathbb P}
\\
   \node{E\times F} \arrow{n,l,J}{(Id_E,\pi^{-1})} \arrow[3]{e,r}{
   \mathbb P \, \circ \, (Id_E, \pi^{-1})} \arrow{s,l}{(\pi,Id_F)}
   \node[3]{\mathbb R_+\cup \{0\}}
\\
 \node{F\times F} \arrow{neee,b,--}{P^*}
\end{diagram}
\end{equation}
where $Id_E:E\to E$ is the identity map on $E$, $\mathfrak P (E)$
is the power set of $E$ and $\mathbb{P} (e,A)=\sum_{e_i\in A}
P((e,e_i))$.
\end{itemize}
When the cardinality of $F$ is strictly less than the cardinality
of $E$, we have reduced some information: the subsets
$\pi^{-1}(f),f\in F$ of $E$ act in the same way for the target
problem.

The proof of the optimal solution's existence was therefore based
on the fact that the set of compatible projections $\pi$ has a
minimal majorant property. More precisely, if $\pi$ is a
projection from $E$ to another set $F$, let $R_\pi$ be the
equivalence relationship on $E$ defined by $e_1\,R_\pi\,e_2 \iff
\pi(e_1)=\pi(e_2)$. If we define
$\widehat{E}:=\{R_\pi\colon\pi\text{ satisfies
}\ref{eq:diagram}\}$, in \cite{AM} it was proved that
$$
{\bigcap \{ R\colon R_\pi\subseteq R,\forall
R_\pi\in\widehat{E}\}} \in \widehat{E} \,.
$$
Unfortunately, finding a nontrivial $R_\pi\in\widehat{E}$ is not a
local search. In fact, we may have $ P(e_1,e_3) \neq P(e_2,e_3) $
but $\mathbb P(e_1,\{e_3,e_4\}) =\mathbb P(e_2,\{e_3,e_4\}) $,
which means that $e_1\,R_\pi\,e_2$ may be found if \emph{we know}
that $e_3\,R_\pi\,e_4$. Moreover, it is not difficult to build
examples where the only nontrivial element of $\widehat{E}$
corresponds to the optimal nontrivial projection. Therefore,
searching for a compressing map $\pi$ appears as a non--polynomial
search, in the sense that we have to look at the whole set
$\widetilde{E}$ of equivalent relations on $E$. In fact, finding a
reducing map means to find $R\in \widetilde{E}$ s.t.
\begin{itemize}
\item $\forall e_i\in T$ $e_i\,R\,e_j\iff e_j\in T$; \item
$\forall \{e_i ,e_j, e_k\}\subseteq E \colon e_i R e_j $,
$\sum_{e_l\,R\,e_k} P(e_i,e_l) = \sum_{e_l\,R\,e_k} P(e_j,e_l) $.
\end{itemize}
The problem here is to find a polynomial algorithm for reaching
the optimal projection of the given Markov chain $(E,P)$ which
preserves probabilities of reaching the target set $T$. Moreover,
we extend this method to multi--target problems
$\mathbf{T}=\{T_1,\ldots,T_k\}$.

\section{The Target Algorithm}
As in \cite{AM}, we act on the set of equivalence relations on a
set, but we will focus our attention on $F$ instead of on $E$. For
any (finite) set $A$, we denote by $\widetilde{A}$ be the set of
all equivalence relations on $A$. Moreover, given an equivalence
relation $R\in\widetilde{A}$, we denote by $A/R$ the quotient set
of $A$ by $R$. We introduce a partial order $\vDash$ on
$\widetilde{A}$. Let $R,S\in\widetilde{A}$. We say that $R\vDash
S$ if $a_1\,R\,a_2$ implies $a_1\,S\,a_2$ (if you think $A$ as the
set of all men and $R$ is ``belonging to the same state'' while
$S$ is ``belonging to the same continent'', then $R\vDash S$). The
relation $\vDash$ is just set-theoretic inclusion between
equivalence relations, since any relation is a subset of $A\times
A$. We denote by $|A|$ the cardinality of a set $A$. We state the
following trivial lemma without proof.
\begin{lemma}\label{lem:eq_mon}
Let $A$ be a set. $|\cdot|$ is monotone with respect to $\vDash$
in $\widetilde{A}$, i.e.\
\begin{subequations}\label{eq:mon_num}
\begin{gather}\label{eq:mon_numa}
\forall R,S\in\widetilde{A} ,\qquad R\vDash S \Longrightarrow |A/R|\geq |A/S| \,.\\
\intertext{Moreover, if $|A|<\infty$, $|\cdot|$ is strictly
monotone:} |A/R|= |A/S|, R\vDash S \Longrightarrow
R=S\,.\label{eq:mon_numb}
\end{gather}
\end{subequations}
\end{lemma}

Let $(E,T,P)$ be a triple, as above and let $\pi:E\to F$ be the
optimal projection, (for existence and uniqueness, see \cite{AM}).
The map $\pi$ is characterized by the equivalence relationship
$F_\pi$ on $E$: $e_1\,F_\pi\,e_2 \iff \pi(e_1)=\pi(e_2)$.

Let $\widetilde{F}_t$ be the set of all equivalence relations on
$F$ such that the target state $t\in F$ is left ``alone'': i.e.\
$R\in\widetilde{F}_t$ if $t\,R\,f \iff f=t$.

Note that $\widetilde{F}_t \hookrightarrow \widetilde{E}$; more
precisely, since $E\mathop{\longrightarrow}\limits ^{\pi} F$, we
have:
\begin{equation*}
\widetilde{F}_t \mathop{\hookrightarrow}\limits^{j} \mathfrak
P(F\times F) \mathop{\longrightarrow}\limits^{(\pi,\pi)^{-1}}
\mathfrak P(E\times E)
\end{equation*}
It is obvious that $(\pi,\pi)^{-1}\circ j:\widetilde{F}_t\to
\mathfrak P(E\times E)$ defines an equivalence relationship on
$E$. With this inclusion in mind, we can state that
$\widetilde{F}_t \subseteq \widetilde{E}$:
\begin{equation}\label{eq:uniqness}
\widetilde{F}_t \longleftrightarrow \{R\in\widetilde{E}\colon
F_\pi\vDash R\}\,,
\end{equation}
and hence we refer to $\widetilde{F}_t$ both as a class of
equivalence relations on $F$ and on $ E$. The uniqueness of the
optimal solution in \cite{AM} states that \eqref{eq:uniqness} is
well--posed.

We call $I_F$ the identity relationship on $F$:
$$f_1\,I_F\,f_2 \iff f_1=f_2$$i.e.\ $I_F$ is just $F_\pi$ on $\widetilde{F}_t$, and let
$M_E$ be maximal relationship on $\widetilde{F}_t$:
$$
e_1\,M_E\,e_2 \iff \{e_1,e_2\}\subseteq T \text{ or }
\{e_1,e_2\}\subseteq (E\setminus T)\,.
$$
Clearly, $M_E\in\widetilde{F}_t$ and $I_F\vDash R\vDash M_E$,
$\forall R\in\widetilde{F}_t$ (i.e., $I_F$ and $M_E$ are the
minimal and maximal relationship on $\widetilde{F}_t$). Note that
we can compute $M_E$ without knowing $F$.
\medskip

We build now a monotone operator $\mathcal{P}$ on $\widetilde{E}$
(the algorithm's
idea will be to reach $I_F$ --unknown-- starting from $M_E$ --known--).\\
Let $\mathcal{P}:\widetilde{E}\to\widetilde{E}$ so defined:\\
for any $R\in\widetilde{E}$, let $r_1,\ldots,r_N$ be the classes
of equivalence of $E$ induced by $R$. Define
\begin{gather*}
e_1 \mathcal{P}_{r_i} e_2 \iff \mathbb{P}(e_1,r_i)=\mathbb{P}(e_2,r_i) \\
\mathcal{P}(R) = \bigcap_{i=1,\ldots,N} \mathcal{P}_{r_i} \cap
R\,.
\end{gather*}
Now, we focus our attention on the action of $\mathcal{P}$ on
$\widetilde{F}_t$. First, we prove that
$\mathcal{P}:\widetilde{F}_t\to \widetilde{F}_t$ and then we will
show that the unique fixed point of
$\mathcal{P}:\widetilde{F}_t\to \widetilde{F}_t$ is $I_F$.
\begin{lemma}\label{lem:mon}
$\mathcal{P}_{|_{\widetilde{F}_t}}:\widetilde{F}_t\to\widetilde{F}_t$
is a $\vDash$--monotone operator on $\widetilde{F}_t$.
\end{lemma}
\begin{proof}
Let $R\in\widetilde{F}_t\hookrightarrow\widetilde{E}$. Since every
$r_i$ is a subset of $F$, the existence of $P^*$ in
\eqref{eq:diagram} ensures that
$\mathcal{P}_{r_i}\in\widetilde{F}_t$. Therefore,
$\mathcal{P}(R)\in\widetilde{F}_t$. Since $\mathcal{P}(R)
\subseteq R$, it is a monotone operator.
\end{proof}
\begin{theorem}\label{thm:monotone}
$I_F$ is the unique fixed point of
$\mathcal{P}:\widetilde{F}_t\to\widetilde{F}_t$.
\end{theorem}
\begin{proof}
$F_\pi$ is trivially a fixed point for
$\mathcal{P}:\widetilde{E}\to\widetilde{E}$ by \eqref{eq:diagram}
and then $I_F$ is a fixed point for
$\mathcal{P}:\widetilde{F}_t\to\widetilde{F}_t$.

Now, let $R\in\widetilde{F}_t$ s.t.\ $R=\mathcal{P}(R)$. Define
the canonical map $\pi_R:F\to F/R$. We have
\begin{itemize}
    \item $(\pi_R\circ\pi) (E) = F/R$;
    \item $(\pi_R\circ\pi) (T) =\pi_R(\pi(T))=\pi_R(t)=t$, \\
$(\pi_R\circ\pi)^{-1} (t) =\pi^{-1}(\pi_R^{-1}(t))=\pi^{-1}(t)=T$;
    \item $R\subseteq \mathcal{P}_{r_i}$ $\forall i$, and hence
the following diagram commutes:
\begin{equation*}
\begin{diagram}
   \node{E\times \mathfrak P(E)} \arrow{seee,l}{
   \mathbb P}
\\
   \node{E\times F/R} \arrow{n,l,J}{(Id_E,(\pi_R\circ\pi)^{-1})} \arrow[3]{e,r}{
   \mathbb P \, \circ \, (Id_E, (\pi_R\circ\pi)^{-1})} \arrow{s,l}{(\pi_R\circ\pi ,I_{F/R})}
   \node[3]{\mathbb R_+\cup \{0\}}
\\
 \node{F/R\times F/R} \arrow{neee,--}
\end{diagram}
\end{equation*}
\end{itemize}
Since $\pi$ is the optimal projection such that \eqref{eq:diagram}
holds, then $F/R\equiv F$, i.e.\ $R=I_F$.
\end{proof}

\begin{corollary}\label{cor:alg}
Let $(E,T,P)$ be given and let $N$ be the cardinality of $F$,
i.e.\ $N=|E/F_\pi|$. Then $\mathcal{P}^{N-2} ( M_E) = F_\pi$,
where $\mathcal{P}^n := (\mathcal{P}\circ\mathcal{P}^{n-1})$ and
$\mathcal{P}^0$ is the identity operator (i.e.\
$\mathcal{P}^0(R)=R$, $\forall R$).
\end{corollary}
\begin{proof}
First, note that $\mathcal{P}^n(M_E)\in\widetilde{F}_t$ $\forall
n$ (by Lemma~\ref{lem:mon}). Therefore, we may consider
$\mathcal{P}^n :\widetilde{F}_t\to\widetilde{F}_t$. We have $I_F
\vDash \mathcal{P}^{n+1}(M_E)\vDash \mathcal{P}^{n}(M_E)\vDash
M_E$ $\forall n$.

\noindent Let $C_n=|E/(\mathcal{P}^n(M_E))|$. We now prove by
induction on $n$ that
\begin{equation}\label{eq:impl}
\mathcal{P}^n(M_E) \neq I_F \quad \Longrightarrow \quad C_n > n+1
\,.
\end{equation}
For $n=0$, $C_0 =2$ (otherwise $E=T$ and the problem is trivial).
For the induction step, if $\mathcal{P}^n(M_E)\neq I_F$, then $C_n
> n+1$. $\mathcal P$ is a monotone operator, then
$\mathcal{P}^{n+1}(M_E)\vDash \mathcal{P}^{n}(M_E)$ and hence
$C_{n+1}\geq C_n$ by \eqref{eq:mon_numa}. Now, if $C_{n+1}= C_n$,
then $\mathcal{P}^{n+1}(M_E) = \mathcal{P}^{n+1}(M_E)$ by
\eqref{eq:mon_numb} which means that $\mathcal{P}^{n}(M_E)=I_F$ by
Theorem~\ref{thm:monotone}. Therefore, \eqref{eq:impl} holds.

If $ \mathcal{P}^{N-3}(M_E) = I_F $, then $ \mathcal{P}^{N-1}(M_E)
= I_F $ by Theorem~\ref{thm:monotone}. As a consequence of
\eqref{eq:impl}, if $ \mathcal{P}^{N-3}(M_E) \neq I_F $, then
$C_{N-2} \geq N = | F |$. Therefore, since
$\mathcal{P}^{N-1}(M_E)\in\widetilde{F}_t$, we have
$\mathcal{P}^{N-2}(M_E)= I_F$, i.e.\ $\mathcal{P}^{N-2}(M_E)=
F_\pi$.
\end{proof}

\begin{remark}
Note that the operator $\mathcal{P}$ may be computed in a
$|E|$--polynomial time. Corollary~\ref{cor:alg} ensures that
$$
\underbrace{\mathcal P \circ \mathcal P \circ \cdots \circ
\mathcal P}_{ \text{at most $|E/F_\pi|-2$ times } (\leq |E|)}
$$
will reach $F$, given any triple $(E,T,P)$. A \textsc{Matlab}
version of such an algorithm for multitarget $\mathbf{T}$ may be
downloaded at \verb"http://www.mat.unimi.it/~aletti"
\end{remark}

\section{Extension to multiple targets and examples}
The previous results and those in \cite{AM} may be extended to
multiple targets problems. More precisely, let ${X}$ be a
stationary Markov chain on a at finite set $E$ and let $\mathbf{T}
=\{T_1,\ldots, T_k$\} be the absorbing disjoint classes of
targets. Our interstest is engaged by the computation of the
probability of reaching $T_i$ by time $\tau$, given the initial
distribution $\mu$ on $E$. If $(\Omega,\mathcal F,Prob)$ is the
underlying probability space, we are accordingly interested in
\begin{equation}\label{eq:Prob}
\big[Prob ( \cup_{m=0}^\tau \{\omega\in \Omega\colon
X_m(\omega)\in T_i\})\big]_{i=1,\ldots k}
\end{equation}
under the assumption that $Prob(\{X_0=e\})=\mu(e)$\,.

The problem is the following: is there a ``minimum'' set $F$ such
that the problem may be projected to a problem on a Markov chain
on $F$, for any initial distribution $\mu$ on $E$?

The answer is trivial, since each target class $T_i$ defines its
equivalence relationship $I_{F_i}$. It is not difficult to show
that the required set $F$ is defined by
$$
F=E/I_F\,, \qquad \text{where }I_F=\bigcap_{i=1,\ldots,k}I_{F_i}
\,.
$$
\begin{definition}
We call \emph{Markov complexity} of the problem $(E,\mathbf{T},P)$
the cardinality of the optimal set $F$.
\end{definition}
\begin{remark}
The condition $ P^{-1}(0,\infty) \cap (T_i\times E)\subseteq
(T_i\times T_i)$ ensures that each $T_i$ is an absorbing state. In
fact this assumption allows to compute \eqref{eq:Prob} by
$P^\tau$. If we are interested in the probability of being in a
target set $T_i$ at time $\tau$, this condition may be dropped,
leaving the compressing problem unchanged.
\end{remark}
\medskip

\noindent We start here by showing some ``irreducible'' classical
problems.
\begin{example}[Negative Binomial Distribution]
Repeate independently a game with probability $p$ of winning until
you win $n$ games.
\end{example}
Let $S_n = \sum_{i=1}^n Y_i$, where $\{Y_i,i\in\mathbb N\}$ is a
sequence of i.i.d.\ bernoulli random variable with
$Prob(\{Y_i=1\})=1-Prob(\{Y_i=0\})=p$. Our interstest is engaged
by the computation of the probability of reaching $n$ starting
from $0$. Let $E=\{0,1,\ldots,n\}$ be the set of levels we have
reached. We have
\begin{equation*}
\begin{array}{c|cccccc}
   & 0 & 1 & 2 & \ldots & n-1 & n=T \\
   \hline \\[-.4cm]
  0 & (1-p) & p & 0 & \ldots & 0 & 0\\
  1 & 0 & (1-p) & p & \ddots & 0 & 0 \\
  2 & 0 & 0 & (1-p) & \ddots & 0 & 0 \\
  \vdots & \vdots & \vdots & \vdots & \ddots & \ddots & \vdots \\
  n-1 & 0 & 0 & 0 &  \hdots &(1-p) & p \\
  n=T & 0 & 0 & 0 & \hdots & 0 & 1 \\
\end{array}
 \!\!\!
\begin{array}{l}
\begin{array}{l@{}}
\phantom{}
\end{array}
\\
\left.
\begin{array}{l@{}}
\phantom{} \\
\phantom{} \\
\phantom{} \\
\phantom{} \\
\phantom{} \\
\phantom{} \\
\phantom{}
\end{array}
\right\} =:{P}
\end{array}
\end{equation*}
Since the length of the minimum path for reaching the target state
$n$ from different states is different, the problem is irreducible
by \cite[Proposition~31]{AM}. Its Markov complexity is $n+1$.

\begin{example}[Consecutive winning]
Repeate independently a game with probability $p$ of winning until
you win $n$ consecutive games.
\end{example}
The problem is similar to the previous one, where
\begin{equation*}
\begin{array}{c|cccccc}
   & 0 & 1 & 2 & \ldots & n-1 & n=T \\
   \hline \\[-.4cm]
  0 & (1-p) & p & 0 & \ldots & 0 & 0\\
  1 & (1-p) & 0 & p & \ddots & 0 & 0 \\
  2 & (1-p) & 0 & 0 & \ddots & 0 & 0 \\
  \vdots & \vdots & \vdots & \vdots & \ddots & \ddots & \vdots \\
  n-1 & (1-p) & 0 & 0 & \hdots & 0 & p \\
  n=T & 0 & 0 & 0 & \hdots & 0 & 1 \\
\end{array}
 \!\!\!
\begin{array}{l}
\begin{array}{l@{}}
\phantom{}
\end{array}
\\
\left.
\begin{array}{l@{}}
\phantom{} \\
\phantom{} \\
\phantom{} \\
\phantom{} \\
\phantom{} \\
\phantom{} \\
\phantom{}
\end{array}
\right\} =:{P}
\end{array}
\end{equation*}
The problem is again irreducible by \cite[Proposition~31]{AM}. Its
Markov complexity is $n+1$.

\begin{example}[Gambler's ruin]
Let two players each have a finite number of pennies (say, $n_1$
for player one and $n_2$ for player two). Now, flip one of the
pennies (from either player), with the first player having $p$
probability of winning, and transfer a penny from the loser to the
winner. Now repeat the process until one player has all the
pennies.
\end{example}
Let $S_n = \sum_{i=1}^n (2Y_i-1)$, where $\{Y_i,i\in\mathbb N\}$
is a sequence of i.i.d.\ bernoulli random variable with
$Prob(\{Y_i=1\})=1-Prob(\{Y_i=0\})=p$. Our interstest is engaged
by the computation of the probability of reaching $T_1=n_2$ or
$T_2=-n_1$ (multiple target) starting from $0$. Let
$E=\{-n_2,\ldots,-1,0,1,\ldots,n_1\}$ be the set of levels we have
reached. We have
\begin{footnotesize}
\begin{equation*}
\begin{array}{c|ccccccccc}
   & -n_1=T_2 & -n_1+1 & \ldots & -1 & 0 & 1 & \ldots & n_2-1 & n_2=T_1 \\
   \hline \\[-.2cm]
  -n_1=T_2 & 1 & 0 & \ldots & 0 & 0 & 0 & \ldots & 0 & 0\\
  -n_1+1 & (1-p) & 0 & \ddots & 0 & 0 & 0 & \ldots & 0 & 0\\
  \vdots & \vdots & \ddots & \ddots & \ddots & \vdots & \vdots & \vdots & \vdots & \vdots \\
  -1 & 0 & 0 & \ddots & 0 & p & 0 & \ldots & 0 & 0\\
  0 & 0 & 0 & \vdots & (1-p) & 0 & p & \ldots & 0 & 0\\
  1 & 0 & 0 & \vdots & 0 & (1-p) & 0 & \ddots & 0 & 0\\
  \vdots & \vdots & \vdots & \vdots & \vdots & \vdots & \ddots & \ddots & \ddots & \vdots \\
  n_2-1 & 0 & 0 & \vdots & 0 & 0 & 0 & \ldots & 0 & p\\
  n_2=T_1 & 0 & 0 & \ldots & 0 & 0 & 0 & \ldots & 0 & 1\\
\end{array}
\end{equation*}
\end{footnotesize}
This problem is clearly irreducible, since it is for $T_1$ (for
example). The problem may be reduced if and only if we are
interesting in the time of stopping (without knowing who wins,
i.e.\ $\mathbf{T}=T_1\cup T_2$) and $p=1/2$. In this case, the
relevant information is the distance from the nearest border and
hence the problem may be half--reduced.
\medskip

\noindent The following classical problem may be reduced.
\begin{example}[Random walk on a cube]
A particle performs a symmetric random walk on the vertices of a
unit cube, i.e., the eight possible positions of the particle are
$(0, 0, 0)$, $(1, 0, 0)$, $(0, 1, 0)$, $(0, 0, 1)$, $(1, 1,
0)$,\ldots,$(1, 1, 1)$, and from its current position, the
particle has a probability of $1/3$ of moving to each of the $3$
neighboring vertices. This process ends when the particle reaches
$(0, 0, 0)$ or $(1, 1, 1)$.
\end{example}
Let $T_1=(0, 0, 0)$, $T_2=(1, 1, 1)$. The following transiction
matrix
\begin{footnotesize}
\begin{equation*}
\begin{array}{c|cccccccc}
& (0, 0, 0)&(1, 0, 0)&(0, 1, 0)&(0, 0, 1)&(1, 1, 0) &(1, 0, 1) &
(0, 1, 1)&(1, 1, 1)
\\
   \hline \\[-.2cm]
(0, 0, 0) & 1 & 0 & 0 & 0 & 0 & 0 & 0 & 0 \\
(1, 0, 0) & 1/3 & 0 & 0 & 0 & 1/3 & 1/3 & 0 & 0 \\
(0, 1, 0) & 1/3 & 0 & 0 & 0 & 1/3 & 0 & 1/3 & 0 \\
(0, 0, 1) & 1/3 & 0 & 0 & 0 & 0 & 1/3 & 1/3 & 0 \\
(1, 1, 0) & 0 & 1/3 & 1/3 & 0 & 0 & 0 & 0 & 1/3 \\
(1, 0, 1) & 0 & 1/3 & 0 & 1/3 & 0 & 0 & 0 & 1/3 \\
(0, 1, 1) & 0 & 0 & 1/3 & 1/3 & 0 & 0 & 0 & 1/3 \\
(1, 1, 1) & 0 & 0 & 0 & 0 & 0 & 0 & 0 & 1 \\
\end{array}
\end{equation*}
\end{footnotesize}
can be easily reduced on
\begin{equation*}
\begin{array}{c|cccc}
& t_1 & f_1 & f_2 & t_2
\\
   \hline \\[-.4cm]
t_1 & 1 & 0 & 0 & 0 \\
f_1 & 1/3 & 0 & 2/3 & 0 \\
f_2 & 0 & 2/3 & 0 & 1/3 \\
t_2 & 0 & 0 & 0 & 1 \\
\end{array}
\end{equation*}
where $t_i=T_i$ and $f_i=\{e=(e_1,e_2,e_3):\sum e_j=i\}$, i.e.\
its Markov complexity is $4$. If we are only interesting in the
time of stopping (i.e.\ $\mathbf{T}=T_1\cup T_2$), the previous
problem may be reduced to a geometrical one (Markov complexity
equal to $2$). Clearly, this results hold also for random walk on
a $d$--dimensional cube.

\begin{example}[Coupon Collector's Problem]
Let $n$ objects $ \{e_1,\ldots,e_n\}$ be picked repeatedly with
probability $p_i$ that object $e_i$ is picked on a given try, with
$ \sum_{i}p_i=1$. Find the earliest time at which all $n$ objects
have been picked at least once.
\end{example}
Let $\Lambda$ be the set of permutations of the $n$ objects. For a
fixed permutation $\lambda=
(e_{\lambda_1},e_{\lambda_2},\ldots,e_{\lambda_n}) \in \Lambda $
we denote by $E^{\lambda^i} =
\{e_{\lambda_1},e_{\lambda_2},\ldots,e_{\lambda_i}\} $ the set of
the first $i$-objects in $\lambda$ (without order!).

Now, let $\mathcal A_\lambda$ be the set of all the paths that
have picked all the $n$ objects with the order given by $\lambda$.
In Pattern--Matching Algorithms framework (see \cite[Section~3 and
Remark~18]{AM}), the \emph{stopping $\lambda$-rule} we consider
here is denoted by
$$
T_\lambda = e_{\lambda_1} \{E^{\lambda^1}\}^\ast e_{\lambda_2}
\{E^{\lambda^2}\}^\ast \cdots e_{\lambda_{n-1}}
\{E^{\lambda^{n-1}}\}^\ast e_{\lambda_{n}} \,,
$$
and it becomes a target state of an enbedded Markov problem on a
graph (see \cite[Section~3]{AM}). The stopping class for the
Coupon Collector's Problem is accordingly $\mathbf{T}=
\cup_{\lambda_\in\Lambda} T_\lambda $.

It is not difficult to show that the general Coupon Collector's
Problem may be embedded into a Markow network of $2^n-1$--nodes
(its general Markov hard complexity), where
$E=\{\mathbf{T},\{E^{\lambda^i}\colon\lambda\in\Lambda,1<i<n\}$,
the transition matrix is given by
$$
P(E^{\lambda^i},E^{\zeta^j})=
\left\{%
\begin{array}{ll}
    \sum_{k\in \lambda^i}p_k, & \text{if $E^{\lambda^i}=E^{\zeta^j}$;} \\
    p_k, & \text{if $j=i+1$ and $E^{\zeta^j}=\{E^{\lambda^i},e_k\}$;} \\
    0, & \hbox{otherwise;} \\
\end{array}%
\right.
$$
$\lambda,\zeta\in\Lambda$ and $Prob(\{X^1=e_k\})=p_k$. Note that
this matrix is not in general reducible.

If some $p_i$ are equal, i.e., when some states act with the same
law with respect to the problem, the set $E$ can be projected into
a minor one. The easiest case (namely, $p_i=1/n$ $\forall i$) is
projected into a $n$--state problem:
\begin{equation*}
\begin{array}{c|cccccc}
   & f_1 & f_2 & f_3 & \ldots & f_{n-1} & T \\
   \hline \\[-.4cm]
  f_1 & 1/n & 1-1/n & 0 & \ldots & 0 & 0\\
  f_2 & 0 & 2/n & 1-2/n & \ddots & 0 & 0 \\
  f_3 & 0 & 0 & 3/n & \ddots & 0 & 0 \\
  \vdots & \vdots & \vdots & \vdots & \ddots & \ddots & \vdots \\
  f_{n-1} & 0 & 0 & 0 &  \hdots &(n-1)/n & 1/n \\
  T & 0 & 0 & 0 & \hdots & 0 & 1 \\
\end{array}
 \!\!\!
\begin{array}{l}
\begin{array}{l@{}}
\phantom{}
\end{array}
\\
\left.
\begin{array}{l@{}}
\phantom{} \\
\phantom{} \\
\phantom{} \\
\phantom{} \\
\phantom{} \\
\phantom{} \\
\phantom{}
\end{array}
\right\} =:{P}
\end{array}
\end{equation*}
with $Prob(\{X^1=f_1\})=1$. Here,
$f_i=\{E^{\lambda^i},\lambda\in\Lambda\}$. The problem is again
irreducible by \cite[Proposition~31]{AM} and its Markov complexity
is $n$. In general, when we have $m$ different values of
$\{p_i,i=1,\ldots,n\}$ (namely, $q_1,\ldots,q_m$), if
$n_m=|k\colon p_k=q_m|$, then the Markov complexity can be easily
proven to be $\prod_{k=1}^m (n_k+1) -1$.



\end{document}